\documentclass[final,1p,times]{elsarticle}

\usepackage{amssymb}
\usepackage{amsthm}
\usepackage{amscd}
\usepackage{bm}
\usepackage{amsmath}
\usepackage{amsfonts}
\usepackage{graphicx}
\newtheorem{theorem}{Theorem}[section]
\newtheorem{remark}[theorem]{Remark}
\newtheorem{example}[theorem]{Example}
\newtheorem{corollary}[theorem]{Corollary}
\newtheorem{lemma}{Lemma}[section]
\newtheorem{proposition}[theorem]{Proposition}
\newtheorem{definition}{Definition}[section]
\usepackage{mathrsfs}
\usepackage{titletoc}
\usepackage{color}

\newcommand{\Hdi}{\mathop{H}\limits^{\rightarrow}}
\newcommand{\Hor}{\mathop{H}\limits^{\leftrightarrow}}

\newcommand{\cmmnt}[1]{}

\usepackage{xeCJK}

\usepackage{float,mathrsfs,bm,multirow,graphics,fancyhdr}

\usepackage{geometry}
\journal{***}

\begin{document}

\begin{frontmatter}

\title{Lin-Lu-Yau Ricci curvature on hypergraphs}

\author[muc,bnu]{Yulu Tian\corref{tian}}
\ead{yltian@muc.edu.cn}
\author[bnu]{Liang Zhao}
\address[muc]{College of Science, Minzu University of China, Beijing 100081, China}
\address[bnu]{School of Mathematical Sciences Key Laboratory of Mathematics and Complex Systems of MOE,\\
	Beijing Normal University, Beijing 100875, China}

\cortext[tian]{Corresponding author.}

\begin{abstract}
	In this paper, we introduce a unified framework for defining Lin--Lu--Yau (LLY) Ricci curvature on both undirected and directed hypergraphs. By establishing upper bounds and monotonicity properties for the centered idleness function $\overline{\kappa}_\alpha$, we justify the well-posedness and compatibility of our definitions.
	Furthermore, we prove Bonnet-Myers-type results for undirected and oriented hypergraphs, which highlights the potential of LLY Ricci curvature in hypergraph analysis, particularly in studying geometric and structural properties. Our results extend the foundational definitions of graph Ricci curvature by Ollivier \cite{Ollivier2009ricci} and Lin-Lu-Yau \cite{Lin2011Ricci} to the hypergraph setting.
\end{abstract}


\end{frontmatter}

\section{Introduction}

The hypergraph concept generalizes traditional graph theory by permitting hyperedges that connect arbitrary sets of vertices. This mathematical framework provides a powerful tool for modeling higher-order interactions, with broad applications across science and engineering, including chemical reaction network \cite{Jost2019hypergraph}, social network \cite{deArruda2020social,Yu2023self}, computer vision \cite{Bai2021multi,Han2024vision}, bioinformatics and medical science \cite{Liu2022multi,Saifuddin2023hygnn}, among others. 

From both theoretical and applied viewpoints, developing a discrete curvature theory for non-smooth discrete structures like graphs and hypergraphs has emerged as a significant research direction. For graphs in particular, this investigation has matured considerably, producing deep theoretical results and practical applications.
In the seminal work \cite{Ollivier2009ricci}, Ollivier introduced a notion of Ricci curvature for metric spaces endowed with Markov chains. Building upon this foundation, Lin, Lu, and Yau \cite{Lin2011Ricci} subsequently adapted Ollivier's definition to establish a discrete Ricci curvature for general undirected graphs, which we hereafter refer to as Lin-Lu-Yau (LLY) Ricci curvature. A recent development in this field is the work of Bai, Huang, Lu and Yau \cite{Bai2021sum}, which introduced a limit-free formulation of LLY curvature, thereby enhancing the mathematical foundation of the concept through star coupling curvature, while simultaneously eliminating both parametric dependencies and limiting procedures. For comprehensive treatments of Ollivier-type Ricci curvature in undirected graphs, we direct readers to \cite{Bauer2012ollivier, Benson2021volume, Devriendt2022discrete, Jost2014ollivier, Lin2010Ricci, Munch2019ollivier} and references therein. The directed graph setting was addressed by Yamada \cite{Yamada2019ricci}, who proposed an extension of LLY curvature to directed graphs, computed explicit examples, and derived various estimates. Further advancing this direction, Ozawa, Sakurai, and Yamada \cite{Ozawa2020geometric} developed an alternative generalization for strongly connected directed graphs through the mean transition probability kernel associated with the Chung Laplacian.

Research on Ricci curvatures for hypergraphs remains relatively underdeveloped compared to its graph counterpart, though the field has seen accelerated progress in recent years. In \cite{Leal2021forman}, Forman curvature was investigated for both undirected and directed hypergraphs, with discussions of its applications to experimental networks. In the study of Ollivier-type curvature for hypergraphs, some research efforts have been directed toward extending Ollivier's original definition to hypergraph contexts. One appropriate approach \cite{Asoodeh2018curvature} is to reduce hypergraphs to graphs through clique expansion, thereby enabling the application of Ollivier's definition on graphs. However, it should be noted that clique expansion is an irreversible process, which inevitably leads to the loss of certain structural information inherent in hypergraphs. Recently, \cite{Eidi2020ollivier} and  \cite{Eidi2020edge} examined how to define Ollivier Ricci curvature on directed hypergraphs and applied it to analyze metabolic networks. Meanwhile, \cite{Coupette2022ollivier} proposed multiple random walk-based definitions of Ollivier curvature for undirected hypergraphs and evaluated their utility across various synthetic and real-world networks.

When considering extensions of LLY curvature to hypergraphs, a series of related results has been established in \cite{Akamatsu2022new, Akamatsu2023weak, Ikeda2021coarse}, all concerning undirected hypergraphs. In \cite{Ikeda2021coarse} and \cite{Akamatsu2023weak}, the $\lambda$ Kantorovich difference $KD_\lambda$ on hypergraphs was defined via the submodular Laplacian, without via a random walk. Then $KD_\lambda$ is employed to establish LLY curvature on hypergraphs. Notably, as $\lambda\rightarrow 0$, $KD_\lambda$ differs from the $L_1$ Wasserstein distance $W_1$ in the definition of LLY curvature \cite{Lin2011Ricci} by an infinitesimal term $o(\lambda)$. Meanwhile, \cite{Akamatsu2022new} introduced a new transport distance $W_h$ dependent on a concave function $h$, subsequently using $W_h$, instead of $W_1$, to define LLY curvature on hypergraphs. 

While these works represent appropriate extensions of LLY curvature to undirected hypergraphs, they exhibit two limitations: (i) the definitions are restricted to undirected cases, and (ii) the resulting curvature formulations are computationally abstract for practical hypergraph analysis. Our work makes two key contributions that address these challenges: (i) we establish a unified framework for LLY curvature that encompasses both undirected and directed hypergraphs, and (ii) our approach, based on hypergraph random walks, aligns with Ollivier's original formulation while offering significantly improved computability. This advancement facilitates practical applications of hypergraph LLY curvature across various domains.


\section{Undirected hypergraphs}\label{Undirected}

Let $\mathcal{H}_{un}=(V,H,\bm{w})$ be a weighted undirected hypergraph, where $V=(x_1,x_2,\cdots,x_N)$ is the vertex set, $H=(h_1,h_2,\cdots,h_M)$ is the hyperedge set, and $\bm{w}=(w_{h_1},w_{h_2},\cdots,w_{h_M})\in\mathbb{R}^M_+$ is the weights of hyperedges. When unambiguous, we shall omit the subscripts and denote vertices, hyperedges, and weights simply by $x$ (or $u, v$), $h$, and $w$ respectively. We denote the cardinalities of $V$ and $H$ by $\vert V\vert=N$ and $\vert H\vert=M$ respectively.

If there exists a hyperedge $h\in H$ containing two distinct vertices $u,v\in V$, we say that $u$ and $v$ are adjacent and denote this by $u\sim v$. For any vertex $x$, let $\Gamma(x)$ denote its neighborhood, i.e., the set of vertices adjacent to $x$,
\begin{equation*}
	\Gamma(x)\coloneqq\{z\in V: x\sim z\}.
\end{equation*}
The degree of a vertex $x\in V$ is defined by
\begin{equation*}
	\operatorname{Deg}(x)\coloneqq\sum_{h:x\in h}w_h.
\end{equation*}
Furthermore, if there is a sequence of hyperedges $\{h_{\{uv\}_1},h_{\{uv\}_2},\cdots,h_{\{uv\}_l}\}$ such that $u\in h_{\{uv\}_1}$, $v\in h_{\{uv\}_l}$ and $h_{\{uv\}_j}\cap h_{\{uv\}_{j+1}}\neq\emptyset$ for $1\leq j\leq l-1$, then vertices $u$ and $v$ are connected by the hyperpath $\gamma=\{h_{\{uv\}_1},h_{\{uv\}_2},\cdots,h_{\{uv\}_l}\}$.
\begin{definition}
	If each pair of two distinct vertices in a hypergraph $\mathcal{H}_{un}=(V,H,\bm{w})$ can be connected by a hyperpath, $\mathcal{H}_{un}$ is a connected hypergraph.
\end{definition}

In this section, unless stated otherwise, all hypergraphs are assumed to be undirected and connected, which allows us to define a distance as follows.
    
\begin{definition}
	The distance between two vertices $u, v\in V$ is
	\begin{equation}\label{un-distance def}
		d(u,v)\coloneqq\inf_\gamma\sum_{h\in\gamma}w_h,
	\end{equation}
	where the infimum is taken over all hyperpaths $\gamma$ connecting $u$ and $v$. Consequently, the diameter of a hypergraph $\mathcal{H}_{un}$ is defined as
	\begin{equation*}
		\operatorname{diam}(\mathcal{H}_{un})\coloneqq\max_{u,v\in V}d(u,v).
	\end{equation*}
Moreover, since a hyperedge can contain more than two vertices, the length of a hyperedge $h=\{x_{h_1},x_{h_2},\cdots,x_{h_k}\}\in H$ is defined as its minimum pairwise distance,
	\begin{equation*}
		L(h)\coloneqq\min_{1\leq i<j\leq k}d(x_{h_i},x_{h_j}).
	\end{equation*}
\end{definition}

\begin{definition}
	Let $\mathcal{P}(V)$ be the set of probability measures over the vertex set $V$ defined by
	\begin{equation*}
		\mathcal{P}(V)\coloneqq\left\{\mu:V\to [0,1] \middle| \sum_{x\in V}\mu(x)=1\right\}.
	\end{equation*}
	Given two probability distributions $\mu, \nu\in\mathcal{P}(V)$, the 1-Wasserstein distance between $\mu$ and $\nu$ is defined as
	\begin{equation*}
		W(\mu,\nu)\coloneqq\inf_{\pi\in\Pi(\mu,\nu)}\sum_{x,y\in V}\pi(x,y)d(x,y),
	\end{equation*}
	where $\Pi(\mu,\nu)$ is the set of probability measures on $V\times V$ that have $\mu$ and $\nu$ as their marginals, i.e., the coupling $\pi\in \Pi$ between $\mu$ and $\nu$	satisfies
	\begin{equation*}
		\left\{
		\begin{aligned}
			&\  \sum_{y\in V}\pi(x,y)=\mu(x),\\
			&\  \sum_{x\in V}\pi(x,y)=\nu(y).
		\end{aligned}
		\right.
	\end{equation*}
\end{definition}

It is straightforward to verify that the distance $d: V\times V \rightarrow \mathbb{R}$ satisfies the metric space axioms (see, e.g., \cite{Brezis2018remarks,Brezis1986harmonic,Montrucchio2019kantorovich}). Consequently, the $1$-Wasserstein distance $W$ inherits these properties, satisfying both symmetry and the triangle inequality.

Now, we are ready to define the Ollivier-Ricci curvature of a hyperedge $h\in H$.
\begin{definition}\label{un-or curvature def}
	For a weighted hypergraph $\mathcal{H}_{un}=(V,H,\bm{w})$,  $\forall h=\{x_{h_1},x_{h_2},\cdots,x_{h_k}\}\in H$ and $\alpha\in[0, 1]$, the Ollivier-Ricci curvature of $h$ is defined as
	\begin{equation}\label{un-centered}
		\overline\kappa_\alpha(h)
		\coloneqq\frac{W_1(h)-W_\alpha(h)}{L(h)}.
	\end{equation}
	and the Ricci curvature for any two vertices $u, v\in V$ is defined as
	\begin{equation*}
		\kappa_\alpha(u,v)\coloneqq1-\frac{W(\mu_u^\alpha,\mu_v^\alpha)}{d(u,v)},
	\end{equation*}
	where $W_\alpha(h)\coloneqq\sum_{1\leq i<j\leq k}W(\mu_{x_{h_i}}^\alpha,\mu_{x_{h_j}}^\alpha)$, the probability measures $\mu_{x_{h_i}}^\alpha$ are defined on $V$ as follows
	\begin{equation*}
		\mu_{x_{h_i}}^\alpha(z)\coloneqq\left\{
		\begin{aligned}
			&\  \alpha, &&\text{if } z=x_{h_i},\\
			&\  (1-\alpha)\sum_{h': x_{h_i},z\in h'}\frac{1}{\vert h'\vert-1}\frac{w_{h'}}{\operatorname{Deg}(x_{h_i})}, &&\text{if } z\in\Gamma(x_{h_i}),\\
			&\  0, &&\text{otherwise},
		\end{aligned} \qquad \forall z\in V.
		\right.
	\end{equation*}
\end{definition}

When the hypergraph degenerates to a classical graph (i.e., every hyperedge $h\in H$ satisfies $|h|=2$), Definition \ref{un-or curvature def} coincides with  the Ollivier-Ricci curvature on graphs. Note that, by Definition \ref{un-or curvature def} and \cite{Ma2024evolution}, it is clear that the probability measures $\mu_{x}^\alpha$ is locally Lipschitz in $\mathbb{R}^m_+$ with respect to $\bm{w}$.

Define
\begin{equation}\label{un-or curvature eq}
	\kappa_\alpha(h)\coloneqq1-\frac{W_\alpha(h)}{L(h)}.
\end{equation}
At idleness one ($\alpha=1$), $\mu_x^1=\delta_x$. Hence
\begin{equation*}\label{un-endpoint}
	W_1(h)=L_{\Sigma}(h)\coloneqq
	\sum_{1\leq i<j\leq k}d(x_{h_i},x_{h_j})
	\quad\text{and}\quad
	\kappa_1(h)=1-\frac{L_{\Sigma}(h)}{L(h)}.
\end{equation*}
In particular, if $k\geq3$, then
\begin{equation*}\label{un-endpoint-negative}
	\kappa_1(h)\leq1-\binom{k}{2}<0.
\end{equation*}
Thus, the quotient $\kappa_\alpha(h)/(1-\alpha)$ used for graph edges cannot define a finite hyperedge curvature when the minimum length is combined with the sum $W_\alpha(h)$. Therefore, we take the endpoint-centered curvature
\begin{equation*}
	\overline\kappa_\alpha(h)=\kappa_\alpha(h)-\kappa_1(h)
\end{equation*}
as the Ollivier-Ricci curvature of a hyperedge.
For any pair of distinct vertices, $\kappa_1(u,v)=0$. Consequently,
\begin{equation}\label{un-centered-pair-relation}
	\overline\kappa_\alpha(h)
	=\frac{1}{L(h)}\sum_{1\leq i<j\leq k}
	d(x_{h_i},x_{h_j})\kappa_\alpha(x_{h_i},x_{h_j}).
\end{equation}
Writing $R(h)\coloneqq L_{\Sigma}(h)/L(h)$, we obtain the following relationship
\begin{equation*}\label{un-centered-comparison}
	R(h)\min_{1\leq i<j\leq k}\kappa_\alpha(x_{h_i},x_{h_j})
	\leq\overline\kappa_\alpha(h)
	\leq R(h)\max_{1\leq i<j\leq k}\kappa_\alpha(x_{h_i},x_{h_j}).
\end{equation*}

We now introduce the notion of 1-Lipschitz functions in the hypergraph setting, followed by a key result from optimal transport theory, which is the Kantorovich-Rubinstein duality theorem.
\begin{definition}
	A function $f: V\to\mathbb{R}$ is $1$-Lipschitz if
	\begin{equation*}
		\vert f(u)-f(v)\vert\leq d(u, v),
	\end{equation*}
	for any $u, v\in V$. We denote the set of $1$-Lipschitz functions on $V$ by $1$-Lip.
\end{definition}

\begin{proposition}
	\emph{(Kantorovich-Rubinstein duality, \cite{Bourne2018ollivier, Villani2021topics})} Let $\mathcal{H}_{un}=(V,H,\bm{w})$ be a weighted hypergraph and  $\mu_{x_i}^\alpha, \mu_{x_j}^\alpha$ be two probability measures on $V$. Then
	\begin{equation*}
		W(\mu_{x_i}^\alpha,\mu_{x_j}^\alpha)=\sup_{f\in 1\text{-Lip}}\sum_{z\in V}f(z)\left(\mu_{x_i}^\alpha(z)-\mu_{x_j}^\alpha(z)\right).
	\end{equation*}
\end{proposition}

\begin{lemma}\label{un-concave}
	For any two vertices $u,v\in V$ and hyperedge $h\in H$, the functions $\kappa_\alpha(u,v)$, $\kappa_\alpha(h)$, and $\overline\kappa_\alpha(h)$ are concave in $\alpha\in[0,1]$.
\end{lemma}

\begin{proof}
	Similar to the method in \cite{Lin2011Ricci}, it is easy to check that $\kappa_\alpha(u,v)$ is concave and $W(\mu_{u}^\alpha,\mu_{v}^\alpha)$ is convex in $\alpha\in[0,1]$, i.e., for $0\leq\alpha<\beta<\gamma\leq1$ and $\lambda=(\gamma-\beta)/(\gamma-\alpha)$, we have
	\begin{equation*}
		W(\mu_{u}^\beta,\mu_{v}^\beta)\leq\lambda W(\mu_{u}^\alpha,\mu_{v}^\alpha)+(1-\lambda)W(\mu_{u}^\gamma,\mu_{v}^\gamma).
	\end{equation*}
	Hence, $W_\alpha(h)=\sum_{1\leq i<j\leq k}W(\mu_{x_{h_i}}^\alpha,\mu_{x_{h_j}}^\alpha)$ is also convex in $\alpha\in[0, 1]$. Consequently, we have
	\begin{equation*}
		\begin{aligned}
			\kappa_\beta(h)&=1-\frac{W_\beta(h)}{L(h)}\\
			&\geq\lambda\left[1-\frac{W_\alpha(h)}{L(h)}\right]+(1-\lambda)\left[1-\frac{W_\gamma(h)}{L(h)}\right]\\
			&=\lambda\kappa_\alpha(h)+(1-\lambda)\kappa_\gamma(h).
		\end{aligned}
	\end{equation*}
	Subtracting the constant $\kappa_1(h)$ shows that $\overline\kappa_\alpha(h)$ is concave as well.
\end{proof}

\begin{lemma}\label{un-curvature bound}
	For any $\alpha\in[0,1]$, $u,v\in V$ and $h\in H$, there hold
	\begin{equation*}
		\kappa_\alpha(u,v)\leq(1-\alpha)\frac{2\max_{h\in H}w_h}{d(u,v)}
	\end{equation*}
	and
	\begin{equation*}
		\overline\kappa_\alpha(h)\leq(1-\alpha)\frac{(k-1)k\max_{h\in H}w_h}{L(h)},
	\end{equation*}
where $k=|h|$.
\end{lemma}

\begin{proof}
	Set $M(x)\coloneqq\max_{h:x\in h}w_h$. Since $d(x,z)\leq w_{h'}$ whenever $x,z\in h'$, the coupling prescribed by the random walk gives
	\begin{equation*}
		\begin{aligned}
			W(\delta_x,\mu_x^\alpha)
			&=(1-\alpha)\sum_{z\sim x}\sum_{h':x,z\in h'}
			\frac{w_{h'}d(x,z)}{(|h'|-1)\operatorname{Deg}(x)}\\
			&\leq(1-\alpha)\frac{\sum_{h':x\in h'}w_{h'}^2}{\operatorname{Deg}(x)}
			\leq(1-\alpha)M(x).
		\end{aligned}
	\end{equation*}
	The triangle inequality and $W(\delta_u,\delta_v)=d(u,v)$ therefore imply
	\begin{equation*}
		W(\mu_u^\alpha,\mu_v^\alpha)
		\geq d(u,v)-(1-\alpha)(M(u)+M(v)).
	\end{equation*}
	Consequently,
	\begin{equation}\label{maxkappa}
		\kappa_\alpha(u,v)
		\leq(1-\alpha)\frac{M(u)+M(v)}{d(u,v)}
		\leq(1-\alpha)\frac{2\max_{h\in H}w_h}{d(u,v)}.
	\end{equation}
	Summing the preceding transport lower bound over all pairs in $h$ yields
	\begin{equation*}
		W_1(h)-W_\alpha(h)
		\leq(1-\alpha)(k-1)\sum_{i=1}^kM(x_{h_i}).
	\end{equation*}
	Together with \eqref{un-centered}, this gives
	\begin{equation*}
		\overline\kappa_\alpha(h)
		\leq(1-\alpha)(k-1)\frac{\sum_{i=1}^kM(x_{h_i})}{L(h)}
		\leq(1-\alpha)\frac{k(k-1)\max_{h\in H}w_h}{L(h)}.
	\end{equation*}
\end{proof}

\begin{remark}\label{un-rough bound}
	Write $L_{\min}(h)=L(h)$ and
	$L_{\max}(h)=\max_{i<j}d(x_{h_i},x_{h_j})$.
	The endpoint centering $\overline\kappa_\alpha(h)=\kappa_\alpha(h)-\kappa_1(h)$ is essential when $L(h)=L_{\min}(h)$ or $L(h)=L_{\max}(h)$. If instead one takes $L(h)=L_{\Sigma}(h)$, then $\kappa_1(h)=0$, and the centered quotient $\overline\kappa_\alpha(h)/(1-\alpha)$ agrees with the usual quotient $\kappa_\alpha(h)/(1-\alpha)$.
\end{remark}

If $f$ is concave on $[0,1]$, then the secant quotient $(f(\alpha)-f(1))/(1-\alpha)$ is nondecreasing in $\alpha\in[0,1)$. Lemma \ref{un-concave}, the identities $\overline\kappa_1(h)=0$ and $\kappa_1(u,v)=0$, and Lemma \ref{un-curvature bound} therefore show that the following finite limits exist:
	\begin{equation}\label{un-limit}
		\kappa(h)=\lim_{\alpha\rightarrow1}\frac{\overline\kappa_\alpha(h)}{1-\alpha}
		\quad\text{and}\quad
		\kappa(u,v)=\lim_{\alpha\rightarrow1}\frac{\kappa_\alpha(u,v)}{1-\alpha}.
	\end{equation}

\begin{definition}\label{un-LLY curvature def}
	For a weighted undirected hypergraph $\mathcal{H}_{un}=(V,H,\bm{w})$, the quantities in \eqref{un-limit} are called the Lin--Lu--Yau (LLY) curvatures of the hyperedge $h$ and of the pair $(u,v)$, respectively. Equivalently, they are the negatives of the left derivatives at $\alpha=1$ of the corresponding idleness functions $\overline\kappa_\alpha(h)$ and $\kappa_\alpha(u,v)$.
\end{definition}

Taking $\alpha\rightarrow1$ in \eqref{un-centered-pair-relation} gives the useful identity
\begin{equation}\label{un-LLY-pair-relation}
	\kappa(h)=\frac{1}{L(h)}\sum_{1\leq i<j\leq k}
	d(x_{h_i},x_{h_j})\kappa(x_{h_i},x_{h_j}).
\end{equation}

\begin{remark}
	The works \cite{Akamatsu2023weak,Ikeda2021coarse} likewise define the LLY curvature of hypergraphs, but their approaches are based on the Kantorovich difference and submodular Laplacian, which differ from our random walk formulation. These highly theoretical generalizations render curvature computations intractable even for simple hypergraphs. In contrast, our Definition \ref{un-LLY curvature def} adheres to Lin-Lu-Yau's original framework \cite{Lin2011Ricci}, offering significant computational advantages. Detailed implementation is provided in the following example.
\end{remark}

\begin{example}
	Consider the hypergraph $\mathcal{H}_4$ depicted in Figure \ref{exampleH4}, with the vertex set $V=\{x_1,x_2,x_3,x_4\}$, the hyperedge set $H=\{h_1=\{x_1,x_2,x_3\},h_2=\{x_1,x_4\}\}$ and the uniform weights $w \equiv 1$. The curvature values for both the hyperedge $h_1$ and its constituent vertex pairs are
	\begin{equation*}
		\begin{aligned}
			&\kappa(x_2,x_3)=\frac{3}{2},\\
			&\kappa(x_1,x_2)=\kappa(x_1,x_3)=\frac{1}{2},\\
			&\kappa(h_1)=\frac{5}{2}.
		\end{aligned}
	\end{equation*}
	Indeed, $L(h_1)=1$ and all three pairwise distances equal one, so the last equality follows from \eqref{un-LLY-pair-relation}.
	
	\begin{figure}[H]
		\centering
		\includegraphics[width=0.25\textwidth]{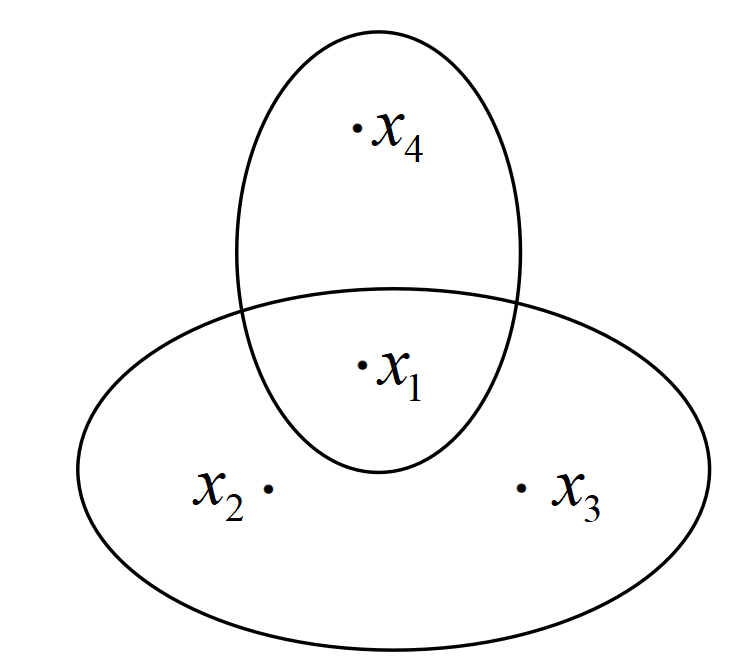}
		\caption{The hypergraph $\mathcal{H}_{4}$}
		\label{exampleH4}
	\end{figure}
\end{example}

When investigating lower bounds of the LLY curvature for arbitrary vertex pairs $(x,y)$, it suffices to consider a special class of pairs defined below, which we term well-transported pairs.
\begin{definition}
	For a hypergraph $\mathcal{H}_{un}=(V,H,\bm{w})$, any pair of vertices $(u,v)$ containing in a hyperedge $h$ is called well-transported if the distance between $u$ and $v$ satisfies $d(u,v)=w_h$.
\end{definition}

\begin{proposition}\label{un-any 2-vertices}
	Suppose $\kappa(u,v)\geq C$ for any well-transported pair $\{u,v\}$, where $C$ is some fixed constant. Then $\kappa(x,y)\geq C$ holds for any two vertices $x,y\in V$.
\end{proposition}

\begin{proof}
	By connectivity of the hypergraph and the finiteness of the hyperedge set, for any two vertices $x,y\in V$, there must exist a sequence of vertices $\{x_{\gamma_1}=x,x_{\gamma_2},\cdots,x_{\gamma_l},x_{\gamma_{l+1}}=y\}$ and an optimal hyperpath $\gamma=\{h_{\gamma_1},h_{\gamma_2},\cdots, h_{\gamma_l}\}$ connecting $x$ and $y$, such that for $j=1,\cdots,l$, the pair $(x_{\gamma_j}, x_{\gamma_{j+1}})$ are well-transported in $h_j$ and $d(x,y)=\sum\limits_{i=1}^l w_{h_{\gamma_i}}$.
	
	By the triangle inequality of $W(\cdot, \cdot)$, for any $\alpha\in[0,1)$, there holds 
	\begin{equation*}
		\begin{aligned}
			\frac{\kappa_\alpha(x,y)}{1-\alpha}&=\frac{1}{1-\alpha}\left[1-\frac{W(\mu_x^\alpha,\mu_y^\alpha)}{d(x,y)}\right]\\
			&\geq\frac{1}{1-\alpha}\left[1-\frac{\sum_{j=1}^lW(\mu_{x_{\gamma_j}}^\alpha,\mu_{x_{\gamma_{j+1}}}^\alpha)}{\sum_{j=1}^lw_{h_{\gamma_j}}}\right]\\
			&=\frac{1}{(1-\alpha)\sum_{j=1}^lw_{h_{\gamma_j}}}\left[\sum_{j=1}^l\left(w_{h_{\gamma_j}}-W(\mu_{x_{\gamma_j}}^\alpha,\mu_{x_{\gamma_{j+1}}}^\alpha)\right)\right]\\
			&=\frac{1}{(1-\alpha)\sum_{j=1}^lw_{h_{\gamma_j}}}\left[\sum_{j=1}^lw_{h_{\gamma_j}}\left(1-\frac{W(\mu_{x_{\gamma_j}}^\alpha,\mu_{x_{\gamma_{j+1}}}^\alpha)}{w_{h_{\gamma_j}}}\right)\right]\\
			&=\frac{\sum_{j=1}^lw_{h_{\gamma_j}}\frac{\kappa_\alpha(x_{\gamma_j},x_{\gamma_{j+1}})}{(1-\alpha)}}{\sum_{j=1}^lw_{h_{\gamma_j}}}.
		\end{aligned}
	\end{equation*}
	Taking limits on both sides as $\alpha\to1$, we conclude that $\kappa(x,y)\geq C$.
\end{proof}

At the end of this section, we prove a Bonnet-Myers type theorem for undirected hypergraphs with positive LLY curvature, which was originally developed for graphs in \cite{Lin2011Ricci,Ollivier2009ricci}.

\begin{theorem}\label{un-positive curvature}
	For any vertices $u,v\in V$, if $\kappa(u,v)>0$, then
	\begin{equation*}
		d(u,v)\leq\frac{2\max_{h\in H}w_h}{\kappa(u,v)}.
	\end{equation*}
	In addition, if there exists $\kappa>0$ such that $\kappa(u,v)\geq\kappa>0$ for any well-transported pair of vertices $(u,v)$, then an upper bound on the hypergraph diameter is
	\begin{equation}\label{un-diam bound}
		\operatorname{diam}(\mathcal{H}_{un})\leq\frac{2\max_{h\in H}w_h}{\kappa}.
	\end{equation}
\end{theorem}

\begin{proof}
	By \eqref{maxkappa} in the proof of Lemma \ref{un-curvature bound}, we have
	\begin{equation*}
		\frac{\kappa_\alpha(u,v)}{1-\alpha}\leq\frac{\max_{h: u\in h}w_h+\max_{h: v\in h}w_h}{d(u,v)}\leq\frac{2\max_{h\in H}w_h}{d(u,v)}.
	\end{equation*}
	Since $\kappa(u,v)>0$, by taking the limit as $\alpha\to1$ we obtain
	\begin{equation*}
		d(u,v)\leq\frac{2\max_{h\in H}w_h}{\kappa(u,v)}.
	\end{equation*}
	If there is a positive lower bound $\kappa$ on the well-transported pair of vertices of the hypergraph, Proposition \ref{un-any 2-vertices} implies that, for any two vertices $x,y\in V$, we have $\kappa(x,y)\geq\kappa$. Therefore, \eqref{un-diam bound} is obtained.
\end{proof}

\section{Directed hypergraphs}\label{Directed}

Let $\mathcal{H}_{di}=(V,\Hdi,\bm{w})$ be a finite and directed hypergraph, where $V=\{x_1,x_2,\cdots,x_N\}$ is the vertex set, $\Hdi=\{h_1,h_2,\cdots,h_M\}$ is the (directed) hyperedge set, and $\bm{w}=(w_{h_1},w_{h_2},\cdots,w_{h_M})\in\mathbb{R}^M_+$ is the weights on hyperedges. We denote the cardinalities of $V$ and $\Hdi$ by $\vert V\vert=N$ and $\vert\Hdi\vert=M$. For a directed hypergraph, each hyperedge $h\in\Hdi$ is an ordered pair $(A_h,B_h)$, where $A_h$ (input) and $B_h$ (output) are nonempty subsets of $V$. If $A_h\cap B_h\neq\emptyset$, we call the hyperedge a (directed) hyperloop.

\begin{definition}
	For two vertices $u$ and $v$, a directed path $\gamma$ from $u$ to $v$ is a sequence of hyperedges $\{h_{\gamma_1},h_{\gamma_2},\cdots,h_{\gamma_l}\}$ such that $u\in A_{h_{\gamma_1}}$, $v\in B_{h_{\gamma_l}}$, and $B_{h_{\gamma_j}}\cap A_{h_{\gamma_{j+1}}}\neq\emptyset$ for $1\leq j\leq l-1$.
\end{definition}


We call a directed hypergraph $\mathcal{H}_{di}=(V,\Hdi,\bm{w})$ weakly connected, if its underlying undirected hypergraph is connected. If for every pair of vertices $u,v\in V$, there exists a directed path from $u$ to $v$, we call the hypergraph strongly connected. In this section, unless stated otherwise, we always suppose that $\mathcal{H}_{di}=(V,\Hdi,\bm{w})$ is directed, loopless (no hyperloops) and strongly connected.

\begin{definition}
	The quasi-distance between two vertices $u, v\in V$ is defined as
	\begin{equation*}
		d(u,v)\coloneqq\inf_\gamma\sum_{h\in\gamma}w_h,
	\end{equation*}
	where the infimum is taken over all directed paths $\gamma$ connecting $u$ and $v$. The diameter of $\mathcal{H}_{di}$ is defined as
	\begin{equation*}
		\operatorname{diam}(\mathcal{H}_{di})\coloneqq\max_{u,v\in V}d(u,v).
	\end{equation*}	
	For a hyperedge $h=(A_h,B_h)\in\Hdi$, with $A_h=\{x^A_1,x^A_2,\cdots,x^A_n\}\xrightarrow{h}B_h=\{y^B_1,y^B_2,\cdots,y^B_m\}$, we define the length of $h$ as
	\begin{equation*}
		L(h)=L(A_h,B_h)\coloneqq\min_{x^A_i\in A_h, y^B_j\in B_h}d(x^A_i,y^B_j).
	\end{equation*}
\end{definition}

It is notable that, although the quasi-distance of two vertices satisfies positivity and the triangle inequality, it need not be symmetric.

\begin{definition}\label{di-neighborhood}
The in-neighborhood of $v\in V$ is
		\begin{equation*}
			\Gamma^{in}(v)\coloneqq\left\{z\in V: \exists h'\in\Hdi, \text{ such that } v\in B_{h'} \text{ and } z\in A_{h'}\right\}.
		\end{equation*}		
The out-neighborhood of $v$ is
		\begin{equation*}
			\Gamma^{out}(v)\coloneqq\left\{z\in V: \exists h'\in\Hdi, \text{ such that } v\in A_{h'}, \text{ and } z\in B_{h'}\right\}.
		\end{equation*}	
Furthermore, for a hyperedge $h=(A_h,B_h)\in\Hdi$, with $A_h=\{x^A_1,x^A_2,\cdots,x^A_n\}$ and $B_h=\{y^B_1,y^B_2,\cdots,y^B_m\}$, the in-neighborhood of $A_h$ is
		\begin{equation*}
			\Gamma^{in}(A_h)\coloneqq\left\{z\in V: \exists h'\in\Hdi, \text{ such that } z\in A_{h'} \text{ and } \exists x^A_i\in B_{h'}\right\}.
		\end{equation*}
The out-neighborhood of $B_h$ is
		\begin{equation*}
			\Gamma^{out}(B_h)\coloneqq\left\{z\in V: \exists h'\in\Hdi, \text{ such that } z\in B_{h'}  \text{ and } \exists y^B_j\in A_{h'}\right\}.
		\end{equation*}
\end{definition}

For a hyperedge $h\in\Hdi$ with $A_h=\{x^A_1,\cdots,x^A_n\}$ and $B_h=\{y^B_1,\cdots,y^B_m\}$, the probability measures $\mu^\alpha_{A_h}$ and $\mu^\alpha_{B_h}$ are defined on $V$ as follows
\begin{equation*}
	\mu^\alpha_{A_h}\coloneqq\sum_{i=1}^{n}\mu^\alpha_{x_i^A},
\end{equation*}
and for any $z\in V$,
\begin{equation*}\label{di-in-measure}
	\mu_{x_i^A}^\alpha(z)\coloneqq\left\{
	\begin{aligned}
		&\  \frac{\alpha}{n}, &&\text{if } z=x_i^A,\\
		&\  (1-\alpha)\sum_{h': x_i^A\in B_{h'},z\in A_{h'}}\frac{1}{n\vert A_{h'}\vert}\frac{w_{h'}}{\sum_{h':x_i^A\in B_{h'}}w_{h'}}, &&\text{if } z\in \Gamma^{in}(x_i^A),\\
		&\  0, &&\text{otherwise}.
	\end{aligned}
	\right.
\end{equation*}
Similarly,
\begin{equation*}
	\mu^\alpha_{B_h}\coloneqq\sum_{j=1}^{m}\mu^\alpha_{y_j^B},
\end{equation*}
and for any $z\in V$,
\begin{equation*}\label{di-out-measure}
	\mu_{y_j^B}^\alpha(z)\coloneqq\left\{
	\begin{aligned}
		&\  \frac{\alpha}{m}, &&\text{if } z=y_j^B,\\
		&\  (1-\alpha)\sum_{h': y_j^B\in A_{h'},z\in B_{h'}}\frac{1}{m\vert B_{h'}\vert}\frac{w_{h'}}{\sum_{h':y_j^B\in A_{h'}}w_{h'}}, &&\text{if } z\in \Gamma^{out}(y_j^B),\\
		&\  0, &&\text{otherwise}.
	\end{aligned}
	\right.
\end{equation*}
The $1$-Wasserstein  quasi-distance between these two measures given by
\begin{equation}\label{di-Wasserstein}
	W(\mu^\alpha_{A_h},\mu^\alpha_{B_h})\coloneqq\inf_{\pi\in\Pi(\mu^\alpha_{A_h},\mu^\alpha_{B_h})}\sum_{u\in V}\sum_{v\in V}\pi(u,v)d(u,v),
\end{equation}
where the coupling $\pi$ between $\mu^\alpha_{A_h}$ and $\mu^\alpha_{B_h}$ satisfies
\begin{equation*}\label{di-coupling}
	\left\{
	\begin{aligned}
		&\  \sum_{v\in V}\pi(u,v)=\sum_{i=1}^{n}\mu^\alpha_{x_i^A}(u),\\
		&\  \sum_{u\in V}\pi(u,v)=\sum_{j=1}^{m}\mu^\alpha_{y_j^B}(v),
	\end{aligned}
	\right.
\end{equation*}
and it is called optimal if it attains the infimum in \eqref{di-Wasserstein}.

One can easily verify that the $1$-Wasserstein quasi-distance inherits positivity and triangle inequality, but it also may not be symmetric. By Definitions \ref{di-neighborhood},  $\forall x_i^A\in A_h$ and $y_j^B\in B_h$, there hold the following identities:
\begin{equation}\label{di-identity}
	\sum_{z\in\Gamma^{in}(x_i^A)}\mu_{x_i^A}^\alpha(z)=\frac{1-\alpha}{n} \quad \text{ and } \quad \sum_{z\in\Gamma^{out}(y_j^B)}\mu_{y_j^B}^\alpha(z)=\frac{1-\alpha}{m}.
\end{equation}Set
\begin{equation*}\label{di-endpoint-measures}
	U_{A_h}\coloneqq\mu_{A_h}^1=\frac1n\sum_{i=1}^n\delta_{x_i^A},
	\quad
	U_{B_h}\coloneqq\mu_{B_h}^1=\frac1m\sum_{j=1}^m\delta_{y_j^B},
\end{equation*}
and let $D_1(h)\coloneqq W(U_{A_h},U_{B_h})$. Every transport from $U_{A_h}$ to $U_{B_h}$ has cost at least $L(h)$, and therefore
\begin{equation}\label{di-endpoint-curvature}
	D_1(h)\geq L(h),
	\quad
	\kappa_1(h)=1-\frac{D_1(h)}{L(h)}\leq0.
\end{equation}
Equality $D_1(h)=L(h)$ holds if and only if there is a coupling of $U_{A_h}$ and $U_{B_h}$ supported on the pairs $(x,y)\in A_h\times B_h$ satisfying $d(x,y)=L(h)$. It holds, in particular, for every loopless unweighted directed hyperedge, since then $d(x,y)=L(h)=1$ for all $x\in A_h$ and $y\in B_h$. In general the inequality in \eqref{di-endpoint-curvature} can be strict, in which case the uncentered quotient $\kappa_\alpha(h)/(1-\alpha)$ tends to $-\infty$,
where
\begin{equation*}
	\kappa_\alpha(h)\coloneqq1-\frac{W(\mu^\alpha_{A_h},\mu^\alpha_{B_h})}{L(h)}.
\end{equation*}

\begin{definition}\label{di-Wasserstein def}
	We define the Ollivier-Ricci curvature of a hyperedge $h\in\Hdi$ with $A_h=\{x^A_1,\cdots,x^A_n\}$ and $B_h=\{y^B_1,\cdots,y^B_m\}$ as follows
	\begin{equation}\label{di-centered}
		\overline\kappa_\alpha(h)
		\coloneqq\kappa_\alpha(h)-\kappa_1(h)
		=\frac{D_1(h)-W(\mu_{A_h}^\alpha,\mu_{B_h}^\alpha)}{L(h)}.
	\end{equation}
\end{definition}


\begin{lemma}\label{di-concave}
	For any hyperedge $h\in\Hdi$, both $\kappa_\alpha(h)$ and $\overline\kappa_\alpha(h)$ are concave in $\alpha\in[0, 1]$.
\end{lemma}

\begin{proof}
	For $0\leq\alpha<\beta<\gamma\leq1$, let $\lambda\coloneqq(\gamma-\beta)/(\gamma-\alpha)$. Suppose that $\pi_\alpha$ and $\pi_\gamma$ are the optimal couplings between $\mu^\alpha_{A_h}, \mu^\alpha_{B_h}$ and $\mu^\gamma_{A_h},\mu^\gamma_{B_h}$ respectively. We have
	\begin{equation*}
		\begin{aligned}
			W(\mu^\alpha_{A_h},\mu^\alpha_{B_h})=\sum_{u,v\in V}\pi_\alpha (u,v)d(u,v), \\ W(\mu^\gamma_{A_h},\mu^\gamma_{B_h})=\sum_{u,v\in V}\pi_\gamma (u,v)d(u,v).
		\end{aligned}
	\end{equation*}
	Define $\pi_\beta\coloneqq\lambda \pi_\alpha+(1-\lambda)\pi_\gamma$. We now claim that $\pi_\beta$ is a coupling between $\mu^\beta_{A_h}$ and $\mu^\beta_{B_h}$. Firstly, we have
	\begin{equation*}
		\begin{aligned}
			\sum_{u\in V}\pi_\beta(u,v)&=\sum_{u\in V}\lambda \pi_\alpha(u,v)+(1-\lambda)\pi_\gamma(u,v)\\
			&=\lambda\mu^\alpha_{B_h}(v)+(1-\lambda)\mu^\gamma_{B_h}(v)\\
			&=\sum_{j=1}^{m}\left(\lambda\mu^\alpha_{y_j^B}(v)+(1-\lambda)\mu^\gamma_{y_j^B}(v)\right).
		\end{aligned}
	\end{equation*}
	To verify the claim, it suffices to establish the identity
	\begin{equation}\label{abg}
		\sum_{j=1}^{m}\lambda\mu^\alpha_{y_j^B}(v)+(1-\lambda)\mu^\gamma_{y_j^B}(v)=\mu^\beta_{B_h}(v),\quad \forall v\in V,
	\end{equation}
	which we demonstrate by examining three distinct cases.
	
	Case (i). $v\in B_h$.
	
	In this case, we use $\left\{y_{j_1}^B,\cdots,y_{j_s}^B\right\}$ to denote the subset (possibly empty) of $B_h$ such that $v\in \Gamma^{out}(y_{j_r}^B)$ for some $1\leq r\leq s$:
	\begin{equation*}
		\left\{y_{j_1}^B,\cdots,y_{j_s}^B\right\}\coloneqq\left\{z\in B_h: \exists h'\in\Hdi, \text{ such that } z\in A_{h'}, v\in B_{h'}\right\}.
	\end{equation*}
	Thus, we have
	\begin{equation*}
		\begin{aligned}
			\sum_{j=1}^{m}\left(\lambda\mu^\alpha_{y_j^B}(v)+(1-\lambda)\mu^\gamma_{y_j^B}(v)\right)&=
			\lambda\frac{\alpha}{m}+(1-\lambda)\frac{\gamma}{m}\\
			&+\sum_{r=1}^{s}(\lambda(1-\alpha)+(1-\lambda)(1-\gamma))\sum_{h': y_{j_r}^B\in A_{h'},v\in B_{h'}}\frac{1}{m\vert B_{h'}\vert}\frac{w_{h'}}{\sum_{h':y_{j_r}^B\in A_{h'}}w_{h'}}\\
			&=\frac{\beta}{m}+\sum_{r=1}^{s}(1-\beta)\sum_{h': y_{j_r}^B\in A_{h'},v\in B_{h'}}\frac{1}{m\vert B_{h'}\vert}\frac{w_{h'}}{\sum_{h':y_{j_r}^B\in A_{h'}}w_{h'}}\\
			&=\sum_{j=1}^{m}\mu^\beta_{y_j^B}(v)\\
			&=\mu^\beta_{B_h}(v).
		\end{aligned}
	\end{equation*}
	Case (ii) $v\in\Gamma^{out}(B_h)$.
	
	In this case, the subset
	\begin{equation*}
		\left\{y_{j_1}^B,\cdots,y_{j_s}^B\right\}\coloneqq\left\{z\in B_h: \exists h'\in\Hdi, \text{ such that } z\in A_{h'}, v\in B_{h'}\right\}
	\end{equation*}
 	of $B_h$ must be nonempty and we obtain
	\begin{equation*}
		\begin{aligned}
			\sum_{j=1}^{m}\left(\lambda\mu^\alpha_{y_j^B}(v)+(1-\lambda)\mu^\gamma_{y_j^B}(v)\right)&=
			\sum_{r=1}^{s}(1-\beta)\sum_{h': y_{j_r}^B\in A_{h'},v\in B_{h'}}\frac{1}{m\vert B_{h'}\vert}\frac{w_{h'}}{\sum_{h':y_{j_r}^B\in A_{h'}}w_{h'}}\\
			&=\sum_{j=1}^{m}\mu^\beta_{y_j^B}(v)\\
			&=\mu^\beta_{B_h}(v).
		\end{aligned}
	\end{equation*}
	Case (iii) $v\notin B_h$ and $v\notin\Gamma^{out}(B_h)$.
	
	In this case, since
	\begin{equation*}
		\mu^\alpha_{y_j^B}(v)=\mu^\beta_{y_j^B}(v)=\mu^\gamma_{y_j^B}(v)=0, \quad \forall 1\leq j\leq m,
	\end{equation*}
	the identity \eqref{abg} obviously holds.
	
	Similarly, we can prove that $\sum\limits_{v\in V}\pi_\beta(u,v)=\mu^\beta_{A_h}(u)$ for any $u\in V$ and the claim is proved. Consequently, we have
	\begin{equation*}
		\begin{aligned}
			W(\mu^\beta_{A_h},\mu^\beta_{B_h})&\leq\sum_{u,v\in V}\pi_\beta(u,v)d(u,v)\\
			&=\lambda\sum_{u,v\in V}\pi_\alpha(u,v)d(u,v)+(1-\lambda)\sum_{u,v\in V}\pi_\gamma(u,v)d(u,v)\\
			&=\lambda W(\mu^\alpha_{A_h},\mu^\alpha_{B_h})+(1-\lambda)W(\mu^\gamma_{A_h},\mu^\gamma_{B_h}),
		\end{aligned}
	\end{equation*}
	which implies that
	\begin{equation*}
		\begin{aligned}
			\kappa_\beta(h)&=1-\frac{W(\mu^\beta_{A_h},\mu^\beta_{B_h})}{L(h)}\\
			&\geq\lambda\left[1-\frac{W(\mu^\alpha_{A_h},\mu^\alpha_{B_h})}{L(h)}\right]+(1-\lambda)\left[1-\frac{W(\mu^\gamma_{A_h},\mu^\gamma_{B_h})}{L(h)}\right]\\
			&=\lambda\kappa_\alpha(h)+(1-\lambda)\kappa_\gamma(h).
		\end{aligned}
	\end{equation*}
	Subtracting the constant $\kappa_1(h)$ proves the assertion for $\overline\kappa_\alpha(h)$.
\end{proof}

\begin{proposition}\label{di-duality}
	 Suppose that $\mu^\alpha_{A_h}$ and $\mu^\alpha_{B_h}$ are probability measures determined by a hyperedge $h=(A_h, B_h)\in \Hdi$. Then we have
	\begin{equation*}
		W(\mu^\alpha_{A_h}, \mu^\alpha_{B_h})\geq\sup_{f\in 1\text{-Lip}}\left(\sum_{i=1}^{n}\sum_{z\in V}f(z)\mu^\alpha_{x_i^A}(z)-\sum_{j=1}^{m}\sum_{z\in V}f(z)\mu^\alpha_{y_j^B}(z)\right).
	\end{equation*}
\end{proposition}

\begin{proof}
	The proof is similar to those in \cite{Eidi2020ollivier,Yamada2019ricci}. For the reader's convenience, we provide the details.	
	For any $1$-Lipschitz function $f$ on $V$, we have
	\begin{equation*}
		\begin{aligned}
			\sum_{u\in V}\sum_{v\in V}\pi(u,v)d(u,v)&\geq\sum_{u\in V}\sum_{v\in V}\left(f(u)-f(v)\right)\pi(u,v)\\
			&=\sum_{u\in V}f(u)\sum_{v\in V}\pi(u,v)-\sum_{v\in V}f(v)\sum_{u\in V}\pi(u,v)\\
			&=\sum_{u\in V}f(u)\mu^\alpha_{A_h}(u)-\sum_{v\in V}f(v)\mu^\alpha_{B_h}(v)\\
			&=\sum_{i=1}^{n}\sum_{z\in V}f(z)\mu^\alpha_{x_i^A}(z)-\sum_{j=1}^{m}\sum_{z\in V}f(z)\mu^\alpha_{y_j^B}(z).
		\end{aligned}
	\end{equation*}
	Since this inequality holds for all 1-Lipschitz functions on hypergraphs and the left-hand side is independent of $f$, we obtain the desired result.
\end{proof}

\begin{lemma}\label{di-curvature bound}
	Let $\nu_{A_h}\coloneqq\mu_{A_h}^0$ and $\nu_{B_h}\coloneqq\mu_{B_h}^0$, and set
	\begin{equation*}
		Q_A(h)\coloneqq W(U_{A_h},\nu_{A_h}),
		\quad
		Q_B(h)\coloneqq W(\nu_{B_h},U_{B_h}).
	\end{equation*}
	For every $\alpha\in[0,1]$,
	\begin{equation*}
		\overline\kappa_\alpha(h)
		\leq\frac{1-\alpha}{L(h)}\bigl(Q_A(h)+Q_B(h)\bigr)
		\leq\frac{2(1-\alpha)\operatorname{diam}(\mathcal H_{di})}{L(h)}.
	\end{equation*}
\end{lemma}

\begin{proof}
	The measures depend affinely on the idleness:
	\begin{equation*}
		\mu_{A_h}^\alpha=\alpha U_{A_h}+(1-\alpha)\nu_{A_h},
		\quad
		\mu_{B_h}^\alpha=\alpha U_{B_h}+(1-\alpha)\nu_{B_h}.
	\end{equation*}
	Mixing an identity coupling with an optimal coupling shows that
	\begin{equation*}
		W(U_{A_h},\mu_{A_h}^\alpha)\leq(1-\alpha)Q_A(h),
		\quad
		W(\mu_{B_h}^\alpha,U_{B_h})\leq(1-\alpha)Q_B(h).
	\end{equation*}
	The triangle inequality for the Wasserstein quasi-distance now gives
	\begin{equation*}
		\begin{aligned}
			D_1(h)
			&=W(U_{A_h},U_{B_h})\\
			&\leq W(U_{A_h},\mu_{A_h}^\alpha)
			+W(\mu_{A_h}^\alpha,\mu_{B_h}^\alpha)
			+W(\mu_{B_h}^\alpha,U_{B_h})\\
			&\leq W(\mu_{A_h}^\alpha,\mu_{B_h}^\alpha)
			+(1-\alpha)\bigl(Q_A(h)+Q_B(h)\bigr).
		\end{aligned}
	\end{equation*}
	Equation \eqref{di-centered} proves the first inequality. Since the hypergraph is finite and strongly connected, every transportation cost is at most $\operatorname{diam}(\mathcal H_{di})$, so $Q_A(h),Q_B(h)\leq\operatorname{diam}(\mathcal H_{di})$.
\end{proof}

By Lemma \ref{di-concave}, the quotient
\begin{equation*}
	g(\alpha)=\frac{\overline\kappa_\alpha(h)}{1-\alpha}
\end{equation*}
is nondecreasing on $[0,1)$. Lemma \ref{di-curvature bound} bounds it from above, so its limit is finite.

\begin{definition}\label{di-LLY curvature def}
	Let $\mathcal{H}_{di}=(V,\Hdi,\bm{w})$ be a weighted directed hypergraph. For any hyperedge $h=(A_h,B_h)\in\Hdi$, the limit 
	\begin{equation}\label{di-limit}
		\kappa(h)=\lim_{\alpha\rightarrow1}
		\frac{\overline{\kappa}_\alpha(h)}{1-\alpha}
	\end{equation}
	is called the Lin--Lu--Yau (LLY) curvature of $h$. It is the negative left derivative of $\kappa_\alpha(h)$ at $\alpha=1$ and agrees with the usual uncentered LLY quotient whenever $D_1(h)=L(h)$.
\end{definition}
	
For $x_i^A\in A_h$ and $y_j^B\in B_h$, let
\begin{equation*}
	\rho_{x_i^A}^\alpha\coloneqq n\mu_{x_i^A}^\alpha,
	\quad
	\sigma_{y_j^B}^\alpha\coloneqq m\mu_{y_j^B}^\alpha.
\end{equation*}
These are probability measures associated with the corresponding directed edge
$e_{ij}=(x_i^A,y_j^B)$. Define its idleness curvature and LLY curvature in the usual way using $W(\rho_{x_i^A}^\alpha,\sigma_{y_j^B}^\alpha)$ and $d(x_i^A,y_j^B)$.

\begin{proposition}\label{di-edge-comparison}
	Suppose that $d(x_i^A,y_j^B)=L(h)$ for every $i,j$. This condition holds for every loopless unweighted directed hyperedge. Then
	\begin{equation*}
		\kappa(h)\geq\min_{1\leq i\leq n,\,1\leq j\leq m}\kappa(e_{ij}).
	\end{equation*}
\end{proposition}

\begin{proof}
	For every $i,j$, let $\pi_{ij}$ be an optimal coupling of
	$\rho_{x_i^A}^\alpha$ and $\sigma_{y_j^B}^\alpha$. The average
	$(nm)^{-1}\sum_{i,j}\pi_{ij}$ couples $\mu_{A_h}^\alpha$ and
	$\mu_{B_h}^\alpha$. Hence, as in \cite[Proposition 2.1]{Eidi2020ollivier},
	\begin{equation*}
		W(\mu_{A_h}^\alpha,\mu_{B_h}^\alpha)
		\leq\frac1{nm}\sum_{i,j}W(\rho_{x_i^A}^\alpha,\sigma_{y_j^B}^\alpha)
		\leq\max_{i,j}W(\rho_{x_i^A}^\alpha,\sigma_{y_j^B}^\alpha).
	\end{equation*}
	The common-distance assumption gives
	$\kappa_\alpha(h)\geq\min_{i,j}\kappa_\alpha(e_{ij})$ and also
	$\kappa_1(h)=0$. Dividing by $1-\alpha$ and taking
	$\alpha\rightarrow1$ proves the proposition.
\end{proof}

On the other hand, LLY curvature of a directed hyperedge cannot always be bounded from above by the supremum of edge curvatures of the corresponding directed graph. We present the following example to illustrate this fact.
	
\begin{figure}[H]
	\centering
	\includegraphics[width=0.4\textwidth]{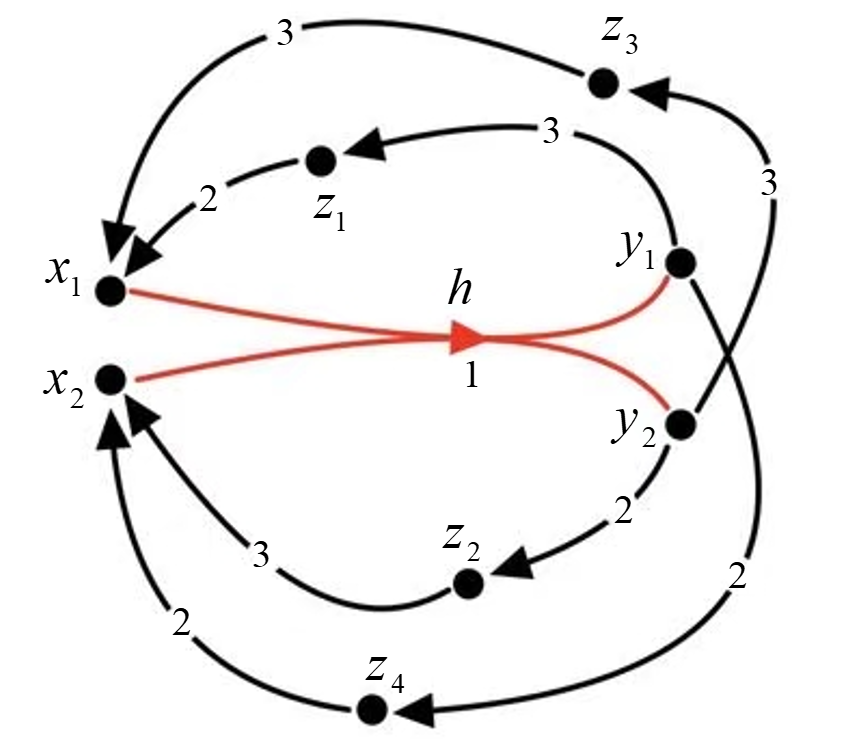}
	\caption{The hypergraph $\mathcal{H}_{xyz}$}
	\label{example2}
\end{figure}
	
\begin{example}
	Consider the hypergraph $\mathcal{H}_{xyz}$ shown in Figure \ref{example2}. Take the hyperedge $h=(A_h,B_h)$ with $A_h=\{x_1,x_2\}$ and $B_h=\{y_1,y_2\}$. Its four corresponding directed edges are $e_{11}=(x_1,y_1)$, $e_{12}=(x_1,y_2)$, $e_{21}=(x_2,y_1)$, and $e_{22}=(x_2,y_2)$. By direct computations we obtain that
	\begin{align*}
		&\kappa(h)=\frac{3}{10}, \\
		&\kappa(e_{11})=-\frac{14}{5}, \quad \kappa(e_{12})=-1, \\
		&\kappa(e_{21})=-\frac{16}{5}, \quad \kappa(e_{22})=-\frac{14}{5}.
	\end{align*}
\end{example}

Since directed hypergraphs may exhibit asymmetric quasi-distances, the induced Wasserstein distance likewise lacks symmetry. When hypergraph distances are symmetric, stronger geometric properties emerge. Henceforth, we focus on oriented hypergraphs with symmetric distance metrics, which has been extensively studied in the spectral theory of hypergraph Laplacians \cite{Andreotti2021Signless, Jost2019hypergraph, Jost2021normalized}.

\begin{definition}\label{or-def}
	A weighted oriented hypergraph is a triple $\mathcal{H}_o=(V,\Hor,\bm{w})$ where:
	\begin{itemize}
		\item $V=(x_1,\cdots,x_N)$ is a finite set of vertices;
		\item $\Hor=(h_1,\cdots,h_M)$ is the set of oriented hyperedges, where each $h\in\Hor$ is an ordered pair $(A_h,B_h)$ with $A_h=\{x^A_1,\cdots,x^A_n\}$ and $B_h=\{y^B_1,\cdots,y^B_m\}$ satisfies $A_h\cap B_h=\emptyset$;
		For every $h=(A_h,B_h)\in\Hor$, its reversal $h^-\coloneqq(B_h,A_h)$ must also belong to $\Hor$;
		\item $\bm{w}: \Hor\to\mathbb{R}_{+}$ is a positive weight function of hyperedges, satisfying the symmetry condition:
		\begin{equation*}
			\bm{w}(h)=\bm{w}(h^-), \quad \forall h\in\Hor.
		\end{equation*}
	\end{itemize}
\end{definition}

The symmetry of an oriented hypergraph sharpens Lemma \ref{di-curvature bound}.

\begin{corollary}\label{di-ez-bounded}
	If $\mathcal{H}_o=(V,\Hor,\bm{w})$ is an oriented hypergraph, then for any $\alpha\in[0,1]$ and $h=(A_h,B_h)\in\Hor$,
	\begin{equation*}
		\overline\kappa_\alpha(h)
		\leq(1-\alpha)\frac{2\max_{e\in\Hor}w_e}{L(h)}.
	\end{equation*}
	Consequently,
	\begin{equation}\label{di-oriented-LLY-bound}
		\kappa(h)\leq\frac{2\max_{e\in\Hor}w_e}{L(h)}.
	\end{equation}
\end{corollary}

\begin{proof}
	In the notation of Lemma \ref{di-curvature bound}, use the random-walk coupling that sends the mass at each $x_i^A$ to its in-neighbors. If $z$ is such an in-neighbor through $e=(A_e,B_e)$, then the reversed hyperedge gives
	$d(x_i^A,z)\leq w_{e^-}=w_e$. Hence
	$Q_A(h)\leq\max_{e\in\Hor}w_e$. The same argument, transporting each out-neighbor of $y_j^B$ back to $y_j^B$, gives
	$Q_B(h)\leq\max_{e\in\Hor}w_e$. Apply Lemma \ref{di-curvature bound} and then let $\alpha\rightarrow1$.
\end{proof}

\begin{remark}\label{di-curvature point}
	(i) For an oriented hypergraph, $\Gamma^{in}(v)=\Gamma^{out}(v)$ for every vertex $v\in V$. We therefore write
	\begin{equation*}
		\Gamma(v)\coloneqq\left\{z\in V: \exists h\in\Hor, v\in A_h, \text{ such that } z\in B_h\right\}.
	\end{equation*}
	Similar to the undirected case, define the Ollivier--Ricci curvature between two distinct vertices of an oriented hypergraph by
	\begin{equation*}
		\kappa_\alpha(u,v)\coloneqq1-\frac{W(\mu^\alpha_{u^{in}},\mu^\alpha_{v^{out}})}{d(u,v)},
	\end{equation*}
	where for any $z\in V$, the probability measures $\mu^\alpha_{u^{in}}$ and $\mu^\alpha_{v^{out}}$ are as follows.
	\begin{equation}\label{di-in-measure-2point}
		\mu_{u^{in}}^\alpha(z)\coloneqq\left\{
		\begin{aligned}
			&\  \alpha, &&\text{if } z=u,\\
			&\  (1-\alpha)\sum_{h': u\in B_{h'},z\in A_{h'}}\frac{1}{\vert A_{h'}\vert}\frac{w_{h'}}{\sum_{h': u\in B_{h'}}w_{h'}}, &&\text{if } z\in \Gamma(u),\\
			&\  0, &&\text{otherwise},
		\end{aligned}
		\right.
	\end{equation}
	\begin{equation}\label{di-out-measure-2point}
		\mu_{v^{out}}^\alpha(z)\coloneqq\left\{
		\begin{aligned}
			&\  \alpha, &&\text{if } z=v,\\
			&\  (1-\alpha)\sum_{h': v\in A_{h'},z\in B_{h'}}\frac{1}{\vert B_{h'}\vert}\frac{w_{h'}}{\sum_{h': v\in A_{h'}}w_{h'}}, &&\text{if } z\in \Gamma(v),\\
			&\  0, &&\text{otherwise}.
		\end{aligned}
		\right.
	\end{equation}
	$W(\mu^\alpha_{u^{in}},\mu^\alpha_{v^{out}})$ is given by
	\begin{equation*}
		W(\mu^\alpha_{u^{in}},\mu^\alpha_{v^{out}})\coloneqq\inf_{\pi\in\Pi(\mu^\alpha_{u^{in}},\mu^\alpha_{v^{out}})}\sum_{a\in V}\sum_{b\in V}\pi(a,b)d(a,b),
	\end{equation*}
	where the coupling $\pi$ between $\mu^\alpha_{u^{in}}$ and $\mu^\alpha_{v^{out}}$ satisfies
	\begin{equation*}
		\left\{
		\begin{aligned}
			&\  \sum_{b\in V}\pi(a,b)=\mu^\alpha_{u^{in}}(a), &&a\in V,\\
			&\  \sum_{a\in V}\pi(a,b)=\mu^\alpha_{v^{out}}(b), &&b\in V.
		\end{aligned}
		\right.
	\end{equation*}
	Furthermore, there holds \cite{Bourne2018ollivier, Villani2021topics}
	\begin{equation*}
		W(\mu^\alpha_{u^{in}},\mu^\alpha_{v^{out}})\geq\sup_{f\in 1\text{-Lip}}\sum_{z\in V}f(z)\left(\mu^\alpha_{u^{in}}(z)-\mu^\alpha_{v^{out}}(z)\right).
	\end{equation*}
	
	(ii) The idleness function $\kappa_\alpha(u,v)$ is concave and satisfies $\kappa_1(u,v)=0$. The same endpoint-stability argument as in Corollary \ref{di-ez-bounded} gives, for $\alpha\in[0,1)$,
	\begin{equation}\label{2point-upperbound}
		\frac{\kappa_\alpha(u,v)}{1-\alpha}\leq\frac{2\max_{h\in\Hor}w_h}{d(u,v)}.
	\end{equation}
	Hence the finite limit
	\begin{equation*}
		\kappa(u,v)\coloneqq\lim_{\alpha\rightarrow1}\frac{\kappa_\alpha(u,v)}{1-\alpha}
	\end{equation*}
	exists and satisfies the same upper bound.
\end{remark}

A Bonnet-Myers type theorem as Theorem  4.1 in \cite{Lin2011Ricci} can be analogously established for oriented hypergraphs.

\begin{theorem}\label{di-BM them}
	Let $\mathcal{H}_o=(V,\Hor,\bm{w})$ be an oriented hypergraph. For any hyperedge $h=(A_h,B_h)\in\Hor$ with $\kappa(h)>0$, there holds
	\begin{equation*}
		L(h)\leq\frac{2\max_{h\in\Hor}w_h}{\kappa(h)}.
	\end{equation*}
	Moreover, if for any hyperedge $h\in \Hor$ and two vertices $u\in A_h, v\in B_h$, there holds $\kappa(u,v)\geq\kappa>0$ for some positive constant $\kappa$, then an upper bound on the hypergraph diameter is
	\begin{equation*}
		\operatorname{diam}(\mathcal{H}_{o})\leq\frac{2\max_{h\in\Hor}w_h}{\kappa}.
	\end{equation*}
\end{theorem}

\begin{proof}
	The first assertion follows immediately from \eqref{di-oriented-LLY-bound}. Taking $\alpha\rightarrow1$ in \eqref{2point-upperbound} gives
	\begin{equation*}
		d(u,v)\leq\frac{2\max_{h\in\Hor}w_h}{\kappa(u,v)}.
	\end{equation*}
	Along a shortest hyperpath, choose consecutive vertices from the tail and head of each traversed oriented hyperedge. These are well-transported pairs, i.e., pairs $(u,v)$ with $u\in A_h$ and $v\in B_h$ such that $d(u,v)=w_h$. The triangle-inequality argument of Proposition \ref{un-any 2-vertices} therefore gives $\kappa(x,y)\geq\kappa$ for every two vertices $x,y\in V$, and the displayed diameter bound follows.
\end{proof}

For simplicity, we always take $\bm{w}\equiv 1$ from now on. As defined in \cite{Andreotti2021Signless,Mulas2021sharp}, the in-degree $\text{deg}_{in}(v)$ of a vertex $v\in V$ is the count of hyperedges containing $v$ as an output, the out-degree $\text{deg}_{out}(v)$ is the count of hyperedges containing $v$ as an input, and the degree of $v$ is 
\begin{equation*}
	\text{deg}(v)\coloneqq\text{deg}_{in}(v)+\text{deg}_{out}(v).
\end{equation*}
For a connected hypergraph, $\operatorname{deg}(v)>0$ for every $v\in V$. Set
$B_H\coloneqq\max_{h\in\Hor}|B_h|$ and
$A_H\coloneqq\max_{h\in\Hor}|A_h|$. Counting incidences gives
	\begin{equation*}
		|\Gamma(v)|\leq B_H\operatorname{deg}_{out}(v)
		\quad\text{and}\quad
		|\Gamma(v)|\leq A_H\operatorname{deg}_{in}(v).
	\end{equation*}

For an oriented hypergraph, we can simplify the probability measures $\mu^\alpha_{A_h}$ and $\mu^\alpha_{B_h}$ on $V$ as follows.
	\begin{equation*}
		\mu^\alpha_{A_h}\coloneqq\sum_{i=1}^{n}\mu^\alpha_{x_i^A},
	\end{equation*}
where for any $z\in V$,
	\begin{equation*}
		\mu_{x_i^A}^\alpha(z)\coloneqq\left\{
		\begin{aligned}
			&\  \frac{\alpha}{n}, &&\text{if } z=x_i^A,\\
			&\  (1-\alpha)\sum_{h': x_i^A\in B_{h'},z\in A_{h'}}\frac{1}{n\vert A_{h'}\vert}\frac{1}{\text{deg}_{in}(x_i^A)}, &&\text{if } z\in \Gamma(x_i^A),\\
			&\  0, &&\text{otherwise}.
		\end{aligned}
		\right.
	\end{equation*}
Similarly,
	\begin{equation*}
		\mu^\alpha_{B_h}\coloneqq\sum_{j=1}^{m}\mu^\alpha_{y_j^B},
	\end{equation*}
where for any $z\in V$,
	\begin{equation*}
		\mu_{y_j^B}^\alpha(z)\coloneqq\left\{
		\begin{aligned}
			&\  \frac{\alpha}{m}, &&\text{if } z=y_j^B,\\
			&\  (1-\alpha)\sum_{h': y_j^B\in A_{h'},z\in B_{h'}}\frac{1}{m\vert B_{h'}\vert}\frac{1}{\text{deg}_{out}(y_j^B)}, &&\text{if } z\in \Gamma(y_j^B),\\
			&\  0, &&\text{otherwise}.
		\end{aligned}
		\right.
	\end{equation*}

For $v\in V$, write
\begin{equation}\label{transition-pv}
	p_v(z)\coloneqq
	\frac{1}{\operatorname{deg}_{out}(v)}
	\sum_{\substack{h:v\in A_h\\ z\in B_h}}\frac1{|B_h|}.
\end{equation}
Thus $p_v$ is the one-step transition probability from $v$. Given distinct vertices $u,v$, put
\begin{equation*}
	\Gamma_u^\varepsilon(v)
	\coloneqq\{z\in\Gamma(v):d(u,z)=d(u,v)+\varepsilon\},
	\quad \varepsilon\in\{-1,0,+1\},
\end{equation*}
and let $p_v^\varepsilon(u)\coloneqq
\sum_{z\in\Gamma_u^\varepsilon(v)}p_v(z)$.

\begin{corollary}\label{di-sym-bounded}
	Let $\mathcal{H}_o=(V,\Hor,\bm{1})$ be an oriented hypergraph and let $u,v\in V$ be distinct. Then
	\begin{equation*}
		\kappa(u,v)\leq
		\frac{1+p_v^{-1}(u)-p_v^{+1}(u)}{d(u,v)}.
	\end{equation*}
\end{corollary}

\begin{proof}
	Put $r=d(u,v)$ and take $\varphi(z)=-d(u,z)$ in the Kantorovich dual lower bound. Since every neighbor of $u$ has distance one from $u$,
	\begin{equation*}
		\begin{aligned}
			W(\mu^\alpha_{u^{in}},\mu^\alpha_{v^{out}})
			&\geq\sum_{z\in V}d(u,z)\mu^\alpha_{v^{out}}(z)
			-\sum_{z\in V}d(u,z)\mu^\alpha_{u^{in}}(z)\\
			&=\alpha r+(1-\alpha)
			\left((r-1)p_v^{-1}(u)+rp_v^0(u)+(r+1)p_v^{+1}(u)\right)
			-(1-\alpha)\\
			&=r-(1-\alpha)
			\left(1+p_v^{-1}(u)-p_v^{+1}(u)\right).
		\end{aligned}
	\end{equation*}
	Divide by $r(1-\alpha)$ and let $\alpha\rightarrow1$.
\end{proof}

Finally, the preceding probability bound yields a volume estimate analogous to
\cite[Theorem 4.3]{Lin2011Ricci}.

\begin{theorem}
	Let $\mathcal{H}_o=(V,\Hor,\bm{1})$ be an unweighted oriented hypergraph. Suppose that $\kappa(u,v)\geq\kappa>0$ whenever $u\in A_h$ and $v\in B_h$ for some $h\in\Hor$. Set
	$\Delta_{out}\coloneqq\max_{v\in V}\operatorname{deg}_{out}(v)$. Then
	\begin{equation*}
		|V|\leq
		1+\sum_{\ell=1}^{\lfloor2/\kappa\rfloor}
		(B_H\Delta_{out})^\ell
		\prod_{i=1}^{\ell-1}\left(1-\frac{i\kappa}{2}\right).
	\end{equation*}
\end{theorem}

\begin{proof}
	Theorem \ref{di-BM them} gives
	$\operatorname{diam}(\mathcal H_o)\leq\lfloor2/\kappa\rfloor$.
	Fix $u\in V$ and write $\Gamma_i(u)=\{z:d(u,z)=i\}$.
	If $v\in\Gamma_i(u)$, then Corollary \ref{di-sym-bounded} and
	$p_v^{-1}(u)+p_v^0(u)+p_v^{+1}(u)=1$ give
	\begin{equation*}
		i\kappa
		\leq i\kappa(u,v)
		\leq1+p_v^{-1}(u)-p_v^{+1}(u)
		\leq2-2p_v^{+1}(u).
	\end{equation*}
	Every $z\in\Gamma_u^{+1}(v)$ occurs in at least one outgoing hyperedge from $v$, so \eqref{transition-pv} implies
	$p_v(z)\geq(B_H\operatorname{deg}_{out}(v))^{-1}$. Therefore
	\begin{equation*}
		|\Gamma_u^{+1}(v)|
		\leq B_H\operatorname{deg}_{out}(v)p_v^{+1}(u)
		\leq B_H\Delta_{out}\left(1-\frac{i\kappa}{2}\right).
	\end{equation*}
	Consequently,
	\begin{equation*}
		|\Gamma_{i+1}(u)|
		\leq|\Gamma_i(u)|B_H\Delta_{out}
		\left(1-\frac{i\kappa}{2}\right).
	\end{equation*}
	Since $|\Gamma_1(u)|\leq B_H\Delta_{out}$, induction yields
	\begin{equation*}
		|\Gamma_\ell(u)|
		\leq(B_H\Delta_{out})^\ell
		\prod_{i=1}^{\ell-1}\left(1-\frac{i\kappa}{2}\right).
	\end{equation*}
	Summing the distance layers proves the result.
\end{proof}

\begin{remark}
	Definition \ref{di-Wasserstein def} and Lemma \ref{di-curvature bound} apply to arbitrary finite strongly connected directed hypergraphs and do not require symmetry. The sharper estimates for oriented hypergraphs use both the symmetry
	\begin{equation}
		d(u,v)=d(v,u), \quad \text{for all } u,v\in V.
	\end{equation}
	and the local in/out structure encoded by Definition \ref{or-def}. They can be extended to any directed hypergraph for which the same symmetry and local transport cost bounds hold. Note that hypergraphs with symmetric distances are not necessarily oriented hypergraphs.
\end{remark}

\begin{remark}\label{un-idleness-linearity}
	Bourne, Cushing, Liu, M{\"u}nch and Peyerimhoff \cite{Bourne2018ollivier} proved that, for adjacent vertices of an unweighted graph equipped with the combinatorial distance and the simple random walk, the idleness function is concave and piecewise linear with at most three linear pieces. Liu and Yang \cite{Liu2026equivalence} proved that, on weighted graphs with general metrics, the idleness function is linear on $[1/2,1]$, and that the threshold $1/2$ is sharp. For unweighted undirected hypergraphs equipped with the combinatorial hyperpath distance and the equal-edges random walk, Xia \cite{Xia2026idleness,Xia2026edge} proved that the idleness function of an adjacent vertex pair is concave and piecewise linear with at most three linear pieces, and is linear on $[1/2,1]$.
	
	The method developed by Liu and Yang \cite{Liu2026equivalence} for weighted graphs can also be adapted to weighted and undirected hypergraphs.  With the distance and probability measures defined in Section \ref{Undirected}, for any $\alpha\in [1/2,1]$, there hold $\kappa_\alpha(h)=(1-\alpha)\kappa(h)$ and $\kappa_\alpha(u,v)=(1-\alpha)\kappa(u,v)$.  Here $h$ denotes an arbitrary hyperedge, and \(u,v\) are an arbitrary pair of vertices that need not lie in the same hyperedge.
\end{remark}




\bibliographystyle{elsarticle-num-names-alpha}

\bibliography{mybib-hypergraph}

@article{Akamatsu2022new,
	title = {A New Transport Distance and Its Associated Ricci Curvature of Hypergraphs},
	author = {Tomoya Akamatsu},
	pages = {90-108},
	volume = {10},
	number = {1},
	journal = {Analysis and Geometry in Metric Spaces},
	year = {2022}
}

@article{Akamatsu2023weak,
	author = {Tomoya Akamatsu},
	title={Weak Kantorovich difference and associated Ricci curvature of hypergraphs},
	journal={arXiv: 2306.14084},
	year={2023}
}

@article{Andreotti2021Signless,
	title={Signless Normalized Laplacian for Hypergraphs},
	author={Andreotti, Eleonora and Mulas, Raffaella},
	journal={arXiv: 2005.14484},
	volume={9},
	year={2021}
}

@INPROCEEDINGS{Asoodeh2018curvature,
	author={Asoodeh, Shahab and Gao, Ting Ran and Evans, James},
	booktitle={2018 IEEE Conference on Decision and Control (CDC)}, 
	title={Curvature of Hypergraphs via Multi-Marginal Optimal Transport}, 
	year={2018},
	pages={1180-1185}
}

@ARTICLE{Bai2021multi,
	author={Bai, Jun Jie and Gong, Biao and Zhao, Yi Ning and Lei, Fu Qiang and Yan, Cheng Gang and Gao, Yue},
	journal={IEEE Transactions on Image Processing}, 
	title={Multi-Scale Representation Learning on Hypergraph for 3D Shape Retrieval and Recognition}, 
	year={2021},
	volume={30},
	pages={5327-5338}
}

@article{Bai2021sum,
	title={On the sum of Ricci-curvatures for weighted graphs},
	author={Bai, Shu Liang and Huang, An and Lu, Lin Yuan and Yau, Shing-Tung},
	journal={Pure Appl. Math. Q.},
	FJOURNAL = {Pure and applied mathematics quarterly},
	volume={17},
	number={5},
	pages={1599--1617},
	year={2021},
	publisher={Int Press Boston, Inc}
}

@article{Bauer2012ollivier,
	title={Ollivier--Ricci curvature and the spectrum of the normalized graph Laplace operator},
	author={Bauer, Frank and Jost, J{\"u}rgen and Liu, Shi Ping},
	journal={Math. Res. Lett.},
	FJOURNAL = {Mathematical Research Letters},
	volume={19},
	number={6},
	pages={1185--1205},
	year={2012},
	publisher={International Press of Boston}
}

@article{Benson2021volume,
	title={Volume growth, curvature, and Buser-type inequalities in graphs},
	author={Benson, Brian and Ralli, Peter and Tetali, Prasad},
	journal={Int. Math. Res. Not. IMRN},
	FJOURNAL = {International Mathematics Research Notices},
	volume={2021},
	number={22},
	pages={17091--17139},
	year={2021},
	publisher={Oxford University Press}
}

@article{Bourne2018ollivier,
	title={Ollivier--Ricci Idleness Functions of Graphs},
	author={Bourne, David P and Cushing, David and Liu, Shi Ping and M{\"u}nch, F and Peyerimhoff, Norbert},
	JOURNAL = {SIAM J. Discrete Math.},
	FJOURNAL = {SIAM Journal on Discrete Mathematics},
	volume={32},
	number={2},
	pages={1408--1424},
	year={2018},
	publisher={SIAM}
}

@article{Brezis2018remarks,
	title={Remarks on the Monge--Kantorovich problem in the discrete setting},
	author={Brezis, Ha{\"\i}m},
	journal={C. R. Acad. Sci. Paris},
	FJOURNAL = {Comptes Rendus. Math{\'e}matique},
	volume={356},
	number={2},
	pages={207--213},
	year={2018}
}

@article{Brezis1986harmonic,
	title={Harmonic maps with defects},
	author={Brezis, Ha{\"\i}m and Coron, Jean-Michel and Lieb, Elliott H},
	journal={Commun. Math. Phys.},
	FJOURNAL = {Communications in mathematical physics},
	volume={107},
	pages={649--705},
	year={1986},
	publisher={Springer}
}

@article{Coupette2022ollivier,
	title={Ollivier-Ricci curvature for hypergraphs: A unified framework},
	author={Coupette, Corinna and Dalleiger, Sebastian and Rieck, Bastian},
	journal={arXiv: 2210.12048},
	year={2022}
}

@article{deArruda2020social,
	title = {Social contagion models on hypergraphs},
	author = {de Arruda, Guilherme Ferraz and Petri, Giovanni and Moreno, Yamir},
	journal = {Phys. Rev. Res.},
	volume = {2},
	issue = {2},
	pages = {023032},
	numpages = {6},
	year = {2020},
	month = {Apr},
	publisher = {American Physical Society}
}

@article{Devriendt2022discrete,
	title={Discrete curvature on graphs from the effective resistance},
	author={Devriendt, Karel and Lambiotte, Renaud},
	journal={J. Phys. Complex.},
	FJOURNAL = {Journal of Physics: Complexity},
	volume={3},
	number={2},
	pages={025008},
	year={2022},
	publisher={IOP Publishing}
}

@article{Eidi2020edge,
	title={Edge-based analysis of networks: curvatures of graphs and hypergraphs},
	author={Eidi, Marzieh and Farzam, Amirhossein and Leal, Wilmer and Samal, Areejit and Jost, J{\"u}rgen},
	journal={Theory Biosci.},
	FJOURNAL = {Theory in biosciences},
	volume={139},
	number={4},
	pages={337--348},
	year={2020},
	publisher={Springer}
}

@article{Eidi2020ollivier,
	title = {Ollivier Ricci curvature of directed hypergraphs},
	author = {Eidi, Marzieh and Jost, J{\"u}rgen},
	JOURNAL = {Sci. Rep.},
	FJOURNAL = {Scientific Reports},
	volume = {10},
	number = {1},
	pages = {12466},
	year = {2020},
}

@INPROCEEDINGS{Han2024vision,
	author={Han, Yan and Wang, Pei Hao and Kundu, Souvik and Ding, Ying and Wang, Zhang Yang},
	booktitle={2023 IEEE/CVF International Conference on Computer Vision (ICCV)}, 
	title={Vision HGNN: An Image is More than a Graph of Nodes}, 
	year={2023},
	pages={19821-19831}
}

@article{Ikeda2021coarse,
	author = {Ikeda, Masahiro and Kitabeppu, Yu and Takai, Yuuki and Ueharam, Takato},
	title={Coarse Ricci curvature of hypergraphs and its generalization},
	journal={arXiv: 2102.00698},
	year={2021}
}

@article{Jost2014ollivier,
	title={Ollivier’s Ricci curvature, local clustering and curvature-dimension inequalities on graphs},
	author={Jost, J{\"u}rgen and Liu, Shiping},
	journal={Discrete Comput. Geom.},
	FJOURNAL = {Discrete \& Computational Geometry},
	volume={51},
	number={2},
	pages={300--322},
	year={2014},
	publisher={Springer}
}

@article{Jost2019hypergraph,
	title={Hypergraph Laplace operators for chemical reaction networks},
	author={Jost, J{\"u}rgen and Mulas, Raffaella},
	journal={Adv. Math.},
	FJOURNAL = {Advances in mathematics},
	volume={351},
	pages={870--896},
	year={2019},
	publisher={Elsevier}
}

@article{Jost2021normalized,
	title={Normalized Laplace operators for hypergraphs with real coefficients},
	author={Jost, J{\"u}rgen and Mulas, Raffaella},
	journal={J. Complex Netw.},
	FJOURNAL={Journal of complex networks},
	volume={9},
	number={1},
	pages={cnab009},
	year={2021},
	publisher={Oxford University Press}
}

@article{Leal2021forman,
	title={Forman--Ricci curvature for hypergraphs},
	author={Leal, Wilmer and Restrepo, Guillermo and Stadler, Peter F and Jost, J{\"u}rgen},
	journal={Adv. Complex Syst.},
	FJOURNAL = {Advances in Complex Systems},
	volume={24},
	number={01},
	pages={2150003},
	year={2021},
	publisher={World Scientific}
}

@article{Lin2011Ricci,
	title={Ricci curvature of graphs},
	author={Lin, Yong and Lu, Lin Yuan and Yau, Shing-Tung},
	journal={Tohoku Math. J.},
	FJOURNAL = {Tohoku Mathematical Journal, Second Series},
	volume={63},
	number={4},
	pages={605--627},
	year={2011},
	publisher={Mathematical Institute, Tohoku University}
}

@article{Lin2010Ricci,
	title={Ricci curvature and eigenvalue estimate on locally finite graphs},
	author={Lin, Yong and Yau, Shing-Tung},
	journal={Math. Res. Lett.},
	FJOURNAL = {Mathematical research letters},
	volume={17},
	number={2},
	pages={343--356},
	year={2010},
	publisher={International Press of Boston}
}

@article{Liu2022multi,
	author = {Liu, Xuan and Song, Cong Zhi and Liu, Shi Chao and Li, Meng Lu and Zhou, Xiong Hui and Zhang, Wen},
	title = {Multi-way relation-enhanced hypergraph representation learning for anti-cancer drug synergy prediction},
	journal = {Bioinformatics},
	volume = {38},
	number = {20},
	pages = {4782-4789},
	year = {2022}
}

@article{Liu2026equivalence,
	title={Equivalence of {Lin--Lu--Yau} curvature and
	{1/2-Ollivier} curvature on weighted graphs},
	author={Liu, Shi Ping and Yang, Yun Yan},
	journal={arXiv: 2608.05939},
	year={2026}
}

@article{Ma2024evolution,
	title={Evolution of weights on a connected finite graph},
	author={Ma, Ji Cheng and Yang, Yun Yan},
	journal={arXiv: 2411.06393},
	year={2024}
}

@article{Montrucchio2019kantorovich,
	title={Kantorovich distance on a weighted graph},
	author={Montrucchio, Luigi and Pistone, Giovanni},
	journal={arXiv: 1905.07547},
	year={2019}
}

@article{Mulas2021sharp,
	title={Sharp bounds for the largest eigenvalue},
	author={Mulas, Raffaella},
	journal={Math. Notes},
	FJOURNAL = {Mathematical notes},
	volume={109},
	number={1},
	pages={102--109},
	year={2021},
	publisher={Springer}
}

@article{Munch2019ollivier,
	title={Ollivier Ricci curvature for general graph Laplacians: heat equation, Laplacian comparison, non-explosion and diameter bounds},
	author={M{\"u}nch, Florentin and Wojciechowski, Rados{\l}aw K},
	journal={Adv. Math.},
	FJOURNAL = {Advances in mathematics},
	volume={356},
	pages={106759},
	year={2019},
	publisher={Elsevier}
}

@article{Ollivier2009ricci,
	title={Ricci curvature of Markov chains on metric spaces},
	author={Ollivier, Yann},
	journal={J. Funct. Anal.},
	FJOURNAL = {Journal of Functional Analysis},
	volume={256},
	number={3},
	pages={810--864},
	year={2009},
	publisher={Elsevier}
}

@article{Ozawa2020geometric,
	title={Geometric and spectral properties of directed graphs under a lower Ricci curvature bound},
	author={Ozawa, Ryunosuke and Sakurai, Yohei and Yamada, Taiki},
	journal={Calc. Var. Partial. Differ. Equ.},
	FJOURNAL = {Calculus of Variations and Partial Differential Equations},
	volume={59},
	pages={1--39},
	year={2020},
	publisher={Springer}
}

@INPROCEEDINGS{Saifuddin2023hygnn,
	author={Saifuddin, Khaled Mohammed and Bumgardner, Briana and Tanvir, Farhan and Akbas, Esra},
	booktitle={2023 IEEE 39th International Conference on Data Engineering (ICDE)}, 
	title={HyGNN: Drug-Drug Interaction Prediction via Hypergraph Neural Network}, 
	year={2023},
	volume={},
	number={},
	pages={1503-1516}
}

@book{Villani2021topics,
	title={Topics in optimal transportation},
	author={Villani, C{\'e}dric},
	volume={58},
	year={2021},
	publisher={American Mathematical Soc.}
}

@article{Xia2026edge,
	title={{Ollivier--Ricci} Idleness Functions and
	Edge-Connectivity of Hypergraphs},
	author={Xia, Qing},
	journal={arXiv: 2608.06029},
	year={2026}
}

@article{Xia2026idleness,
	title={Idleness Functions for {Ollivier--Ricci} Curvature
	on Hypergraphs},
	author={Xia, Qing},
	journal={arXiv: 2608.01970},
	year={2026}
}

@article{Yamada2019ricci,
	title={The Ricci curvature on directed graphs},
	author={Yamada, Taiki},
	journal={J. Korean Math. Soc.},
	FJOURNAL = {Journal of the Korean Mathematical Society},
	volume={56},
	number={1},
	pages={113--125},
	year={2019},
	publisher={The Korean Mathematical Society}
}

@article{Yu2023self,
	title={Self-supervised learning for recommender systems: A survey},
	author={Yu, Jun Liang and Yin, Hong Zhi and Xia, Xin and Chen, Tong and Li, Jun Dong and Huang, Zi},
	journal={IEEE Transactions on Knowledge and Data Engineering},
	volume={36},
	number={1},
	pages={335--355},
	year={2023},
	publisher={IEEE}
}

\end{document}